\documentclass{article}
\usepackage[latin1]{inputenc}
\usepackage{ifthen}
\usepackage{latexsym}
\usepackage{amssymb}
\usepackage{epsf}
\usepackage{graphics}
\usepackage{subfigure}
\usepackage{amsmath,amsfonts,amsthm,amssymb,a4,epsfig,psfrag,enumerate}
\usepackage{listings}
\usepackage[latin1]{inputenc}

\newcommand{\optpar} [2] []
  {\ensuremath{#2\ifempty{#1} {} {\!\left(#1\right)}}}
\newcommand{\optsub} [2] []
  {\ensuremath{#2\ifempty{#1} {} {_{#1}}}}
\newcommand{\optsuper} [2] []
  {\ensuremath{#2\ifempty{#1} {} {^{#1}}}}

\newlength{\boxwd}
\newlength{\boxht}

\newcommand{\Range} {\operatorname{Ran}}

\newcommand{\ifempty} [1] {\ifthenelse{\equal {#1} {} }}

\newcommand\iprod [3] [] {\optsub[\scriptscriptstyle#1]{\ensuremath{\langle#2,#3\rangle}}} 

\newcommand{\dOmega}{\ensuremath{\partial\Omega}}

\newcommand{\Lp} [2] [] {\optpar[#1] {L^{#2}}}                   
\newcommand{\Lone} [1] [] {\Lp[#1]{1}}
\newcommand{\Ltwo} [1] [] {\Lp[#1]{2}}
\newcommand{\Linfty} [1] [] {\Lp[#1]{\infty}}
\newcommand{\Hk} [2] [] {\optpar[#1] {H^{#2}}}
\newcommand{\Hone} [1] [] {\Hk[#1]{1}}

\newcommand{\Htwo} [1] [] {\Hk[#1]{2}}
\newcommand{\Hzerok} [2] [] {\optpar[#1] {H_0^{#2}}}
\newcommand{\Hzeroone} [1] [] {\Hzerok[#1]{1}}

\newcommand{\Rn} [1] [n]{\ensuremath{\mathbb{R}^{#1}}}
\newcommand{\RR} {\ensuremath{\mathbb{R}}}

\newcommand{\lap} {\Delta} 

\newcommand{\seminorm} [1] [\,] {\ensuremath{|{#1}|}}    

\newcommand{\semnrm} [2] [\,]
  {\ensuremath{\seminorm[{#1}]_{\scriptscriptstyle#2}}}
\newcommand{\nrm} [2] [\,]%
  {\ensuremath{|\!|{#1}|\!|_{{\scriptscriptstyle#2}}}}

\newcommand{\Lonenrm} [2] []%
  {\nrm[{#2}] {\Lone[#1]}}
\newcommand{\Ltwonrm} [2] []%
  {\nrm[{#2}] {\Ltwo[#1]}}
\newcommand{\Linftynrm} [2] []%
  {\nrm[{#2}] {\Linfty[#1]}}
\newcommand{\Honenrm} [2] []%
  {\nrm[{#2}] {\Hone[#1]}}
\newcommand{\Htwonrm} [2] []%
  {\nrm[{#2}] {\Htwo[#1]}}

\newcommand{\negnrm} [1] []%
  {\nrm[{#1}] {-1}}
\newcommand{\zeronrm} [1] []%
  {\nrm[{#1}] {0}}
\newcommand{\onenrm} [1] []%
  {\nrm[{#1}] {1}}
\newcommand{\twonrm} [1] []%
  {\nrm[{#1}] {2}}
\newcommand{\inftynrm} [1] []%
  {\nrm[{#1}] {\infty}}

\newcommand{\onesemnrm} [1] []%
  {\semnrm[{#1}] {1}}

\newcommand{\Honesemnrm} [2] []%
  {\semnrm[{#2}] {\Hone[#1]}}
\newcommand{\Htwosemnrm} [2] []%
  {\semnrm[{#2}] {\Htwo[#1]}}

\newtheorem{prop} {Proposition}

\newcommand{\Vp} {\ensuremath{V_{\scriptscriptstyle X}}}%
\newcommand{\Vd} {\ensuremath{V_{\scriptscriptstyle Y}}}%
\newcommand{\Wp} {\ensuremath{W_{\scriptscriptstyle X}}}%
\newcommand{\Wd} {\ensuremath{W_{\scriptscriptstyle Y}}}%
\newcommand{\Pd} {\ensuremath{P_{\scriptscriptstyle Y}}}%
\newcommand{\Pp} {\ensuremath{P_{\scriptscriptstyle X}}}%

\title{A Numerical Algorithm for\\ Ambrosetti-Prodi Type Operators}
\author{Jos\'e Teixeira Cal Neto\footnote{jose.calneto@uniriotec.br, DME, UNIRIO, Rio de Janeiro, Brazil} and Carlos Tomei \footnote{tomei@mat.puc-rio.br, Departamento de Matemática, PUC-Rio, Rio de Janeiro, Brazil} }

\begin{document}
\maketitle

\begin{abstract}
We consider the numerical solution of the equation $- \lap u - f(u) = g$, for the unknown $u$ satisfying Dirichlet conditions in a bounded domain $\Omega$. The nonlinearity $f$  has bounded, continuous derivative.
The algorithm uses the finite element method combined with a global Lyapunov-Schmidt decomposition.
\end{abstract}

\medbreak


\smallbreak

{\noindent\bf MSC-class:} 35B32, 35J91, 65N30.

\section{Introduction}\label{sec:intro}

We consider the partial differential equation
\[ F(u) = - \lap u - f(u) = g, \quad u|_{\partial \Omega} = 0, \]
on domains $\Omega \in \RR^n$, taken to be open, bounded, connected subsets of $\Rn$  with piecewise smooth boundary $\dOmega$, assumed to be at least Lipschitz at all points.
There is a vast literature concerning  the number of solutions for general and positive solutions  for different kinds of nonlinearity $f$ and right-hand side $g$ (to cite a few, \cite{hammerstein:1930}, \cite{dolph:1949}, \cite{ambrosetti:1972}, \cite{manes:1973}, \cite{lazer:1981}, \cite{lupo:1988}, \cite{dancer:1989}, \cite{costa:1992}, \cite{breuermckennaplum:2003}).

Here we assume that the nonlinearity $f:\RR \to \RR$ has a bounded, continuous, derivative, $a \le f'(y) \le b$. We show how a global Lyapunov-Schmidt decomposition introduced by Berger and Podolak \cite{berger:1974} in their proof of the Ambrosetti-Prodi theorem (\cite{ambrosetti:1972}, \cite{manes:1973}) gives rise to a satisfactory solution algorithm using the finite element method. The decomposition was rediscovered by Smiley \cite{smiley:1996}, who realized its potential for numerics: our results advance along these lines.

Write $\lap_D$ for the Dirichlet Laplacian in $\Omega$. The algorithm is especially convenient when the number $d$ of eigenvalues of  $-\lap_D$ in the range of $f'$ is small: the infinite dimensional equation reduces to the inversion of a map from $\RR^d$ to itself.

The subject of semilinear elliptic equations is sufficiently mature that algorithms should stand side by side with theory. The situation may be compared to the study of functions of one variable in a basic calculus course. Some functions, like parabolas, may be handled without substantial computational effort, but  understanding increases with graphs, which are obtained by following a standard procedure.

We do not handle the difficulties and opportunities related to the finite dimensional inversion: a generic solver (as in \cite{smiley:2000} and, for $d=2$, \cite{MST2:1996}) should be replaced by an algorithm which makes use of features inherited by the original map $F$. Here we only deal with examples for which $d=1$ and 2, and  there is some craftsmanship in handling the 2-dimensional example. It is in this step of the PDE solver that delicate issues like nonresonance and lack of properness come up.

The results are part of the PhD thesis of the first author \cite{calneto:thesis}. Complete proofs are presented elsewhere.
The authors are grateful to CAPES, CNPq and Faperj for support.

\section{The basic estimate}\label{sec:basics}
We consider the semilinear elliptic equation presented in the introduction for a nonlinearity $f:\RR \to \RR$ with bounded, continuous derivative.

With these hypotheses, it is not hard to see that $F(u) = -\Delta u - f(u)$ is a $C^1$ map between the Sobolev spaces $H^2_0(\Omega)$ and $L^2(\Omega)=H^0(\Omega)$ and between $H^1_0(\Omega)$ and $H^{-1}(\Omega) \simeq H^1_0(\Omega)$. We concentrate on the second scenario, which is natural for the weak formulation of the problem. Still, the geometric statements below hold in both cases.
To fix notation, set $F: X \to Y$, where $X = H^1_0(\Omega)$ and $Y = H^{-1}(\Omega)$.

The basic estimate is given in Proposition \ref{prop:bound}. Its proof is a simple extension of the argument in \cite{berger:1974}.

Define $\overline{f'(\RR)}=[a,b]$ ($a$ allowed to be $-\infty$) and a larger interval $[\tilde{a},\tilde{b}]\supset[a,b]$. Label the eigenvalues of $-\lap_D$ in non-decreasing order. The \emph{index set} $J$ associated to $[\tilde{a},\tilde{b}]$ is the collection of indices of eigenvalues of $-\lap_D$ in that interval.  The set $J$ is associated to the nonlinearity $f$ if $[\tilde{a},\tilde{b}]=\overline{f'(\RR)}$.
An index set defined this way is \emph{complete}: it contains all indices labeling an eigenvalue in the interval.

Denote the \emph{vertical subspaces} by $\Vp \subset X$ and $\Vd \subset Y$ the spans of the normalized eigenfunctions $\phi_j, j \in J$ in $X$ and $Y$ respectively, with orthogonal complements $\Wp$ and $\Wd$. Let $P$ and $Q$ be the orthogonal projections on $V$ and $W$. Clearly, the dimension of the vertical subspaces equals $|J|$, the cardinality of $J$.
Let $v + \Wp \subset X$ be the \emph{horizontal affine subspace} of vectors $v + w, w \in \Wp$ and consider a \emph{projected restriction} $F_v:v+\Wp\to\Wd$, the restriction of $\Pd F$ to  $v+\Wp$.
\begin{prop} \label{prop:bound}
Let $J$ be the index set associated to the nonlinearity $f$ (or to any interval $[\tilde{a},\tilde{b}]$ containing $\overline{f'(\RR)}$).
Then the derivatives $D F_v : v + \Wp \to \Wd$ are uniformly bounded from below.
More precisely, there exists $C > 0$ such that
\begin{equation}\label{eq:Fvhadamard}
  \forall v \in \Vp \ \forall w\in v+\Wp \ \forall h \in \Wp,\quad\nrm[D F_v(w)h]{Y}\ge C\nrm[h]{X}.
  \end{equation}
All such maps are invertible.
\end{prop}

A direct application of Hadamard globalization theorem (\cite{berger:book}) implies that the projected restrictions are diffeomorphisms, for each $v \in \Vp$.

\section{The underlying picture}\label{sec:geometry}
The geometric implications are very natural. The image under $F$ of each horizontal affine subspace $v + \Wp$ is a \emph{sheet}, i.e., a surface which projects under $\Pd$ diffeomorphically to the horizontal subspace $\Wd$. In particular, every vertical affine subspace $w + \Vd$ intercepts each sheet exactly at a single point. It is not hard to see that the intersection is transversal: tangent spaces of sheet and vertical affine subspace form a direct sum decomposition of $Y$.

A \emph{fiber} is the inverse image of a vertical affine subspace. In a similar fashion, fibers are surfaces of dimension $|J|$ which meet every horizontal affine subspace $v + \Wp$ at a single point --- again, the intersection is transversal. Thus, a vertical subspace parameterizes diffeomorphically each fiber, or, said differently, each fiber has a single point of a given \emph{height}.

Recall a key idea in \cite{berger:1974} and \cite{smiley:2000}. It is clear that $X$ and $Y$ are respectively foliated by fibers and vertical affine subspaces. By definition, all the solutions of $F(u) = g$ must lie in the fiber $\alpha_g = F^{-1}( g + \Vd)$. So, in principle, one might solve the equation by first identifying $\alpha_g \simeq \RR^{|J|}$ and then facing the finite dimensional inversion of $F: \alpha_g \to g + \Vd$.

Horizontal affine subspaces are taken diffeomorphically to sheets, but fibers are not taken diffeomorphically to vertical affine subspaces. In a sense, the nonlinearity of the problem was reduced to a finite dimensional issue.
\begin{figure}[h]
\centering
\mbox{\subfigure{\includegraphics[width=.45\linewidth]{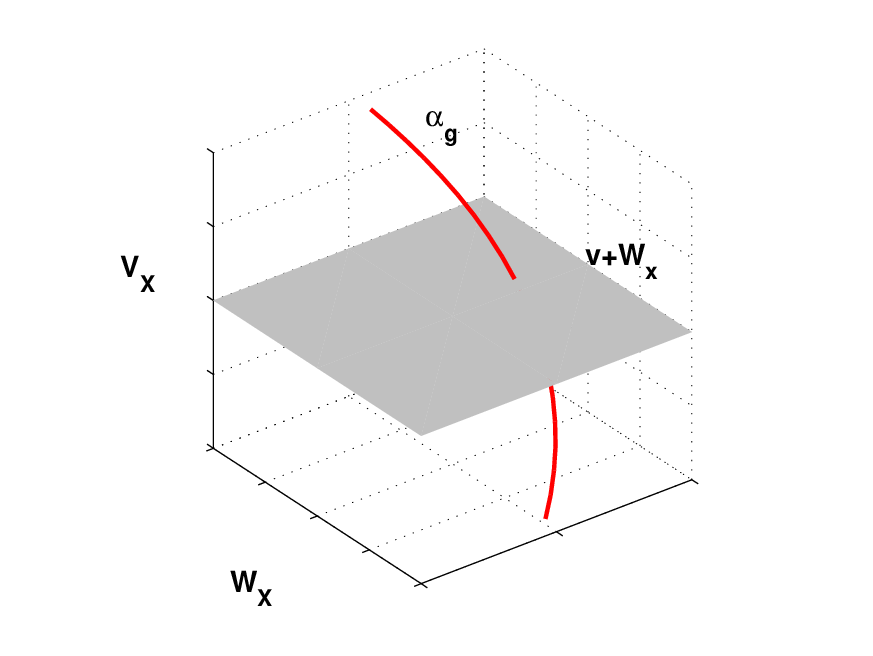}}\quad\quad
\subfigure{\includegraphics[width=.45\linewidth]{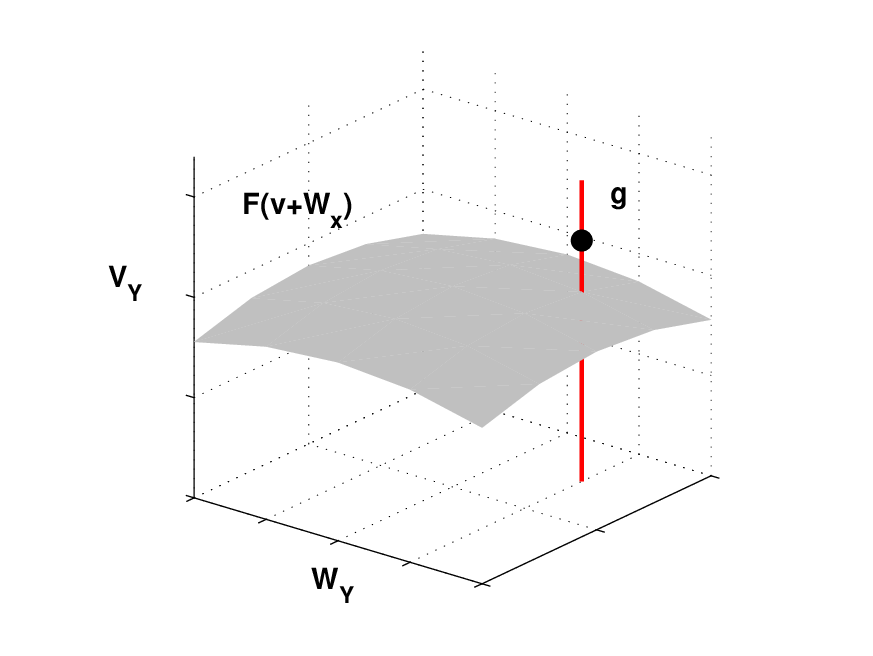}
 }}
\caption{Horizontal affine subspace, fiber; sheet, vertical affine subspace} \label{fig:sheetfib}
\end{figure}

\section{Finding the fiber}\label{sec:findfib}

Recall that each horizontal affine subspace $v + \Wp$ contains exactly one element of each fiber. So, to identify $\alpha_g$, choose $v+\Wp$ and search in it for an element of $\alpha_g$. Said differently, one may think of $F_v: v + \Wp \to \Wd$ as being a diffeomorphism between fibers (represented by points in $v + \Wp$) and vertical affine subspaces (represented by points in $\Wd$).
The situation is ideal for an application of Newton's method: local improvements are performed by linearization of the diffeomorphism.

\begin{figure}[h]
\centering
\mbox{\subfigure{\includegraphics[width=.45\linewidth]{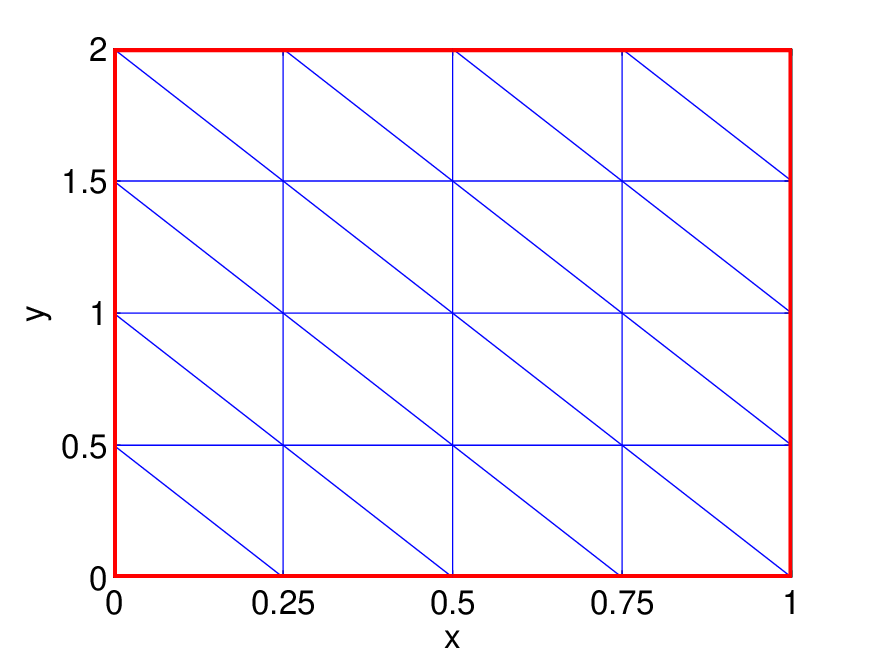}}\quad\quad
\subfigure{\includegraphics[width=.45\linewidth]{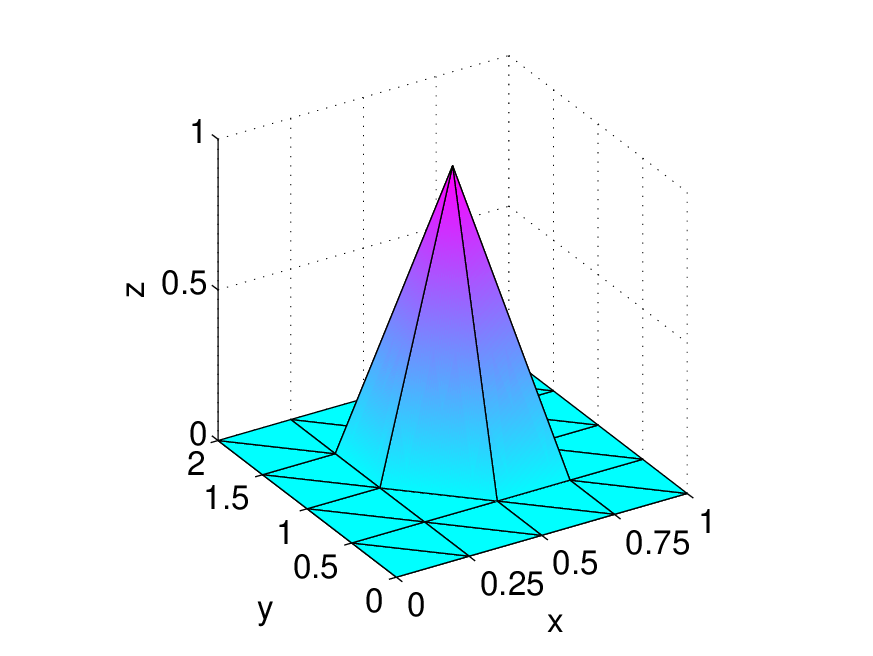}
 }}
\caption{Uniform triangulations of $[0,1]\times[0,2]$ and a nodal function} \label{fig:meshnodalf}
\end{figure}

There is one difficulty, however, related to implementation issues. The functional spaces $X$ and $Y$ give rise to finite dimensional vector spaces generated by \emph{ finite elements}.
We provide some detail: an excellent reference is \cite{ciarlet:book}. First of all, \emph{triangulate} the domain $\Omega$, i.e., split it into disjoint simplices in $\RR^n$. In the examples in Section \ref{sec:eg}, $\Omega$ is the uniformly triangulated rectangle $[0,1]\times[0,2]$.
A \emph{nodal function} is a continuous function that is affine linear on each simplex and has value one at a given vertex and zero at the remaining vertices. These functions form a \emph{nodal basis}, which spans a finite dimensional subspace of \Hzeroone[\Omega]. Figure~\ref{fig:meshnodalf} shows an example of a triangulation of $\Omega$ and one nodal function.

Inner products of nodal functions, both in $\Hone$ and $\Ltwo$ are often zero, a fact which simplifies the numerics associated to the weak formulation of the equation $F(u) = g$.
The vertical subspaces $\Vp$ and $\Vd$ are spanned by eigenfunctions $\phi_j, j \in J$ and are well approximated by a few linear combinations on the nodal basis. On the other hand, obtaining a similar basis for the approximation of the orthogonal subspaces $\Wp$ and $\Wd$ requires much more numerical effort and should be avoided.

To circumvent this problem, extend the Jacobian of $F_v: v + \Wp \to \Wd$ at a point $u$ to an invertible operator $L_u: X \to Y$ which is easy to handle and apply Newton's method to $L_u$ instead.
Setting
  \[
  L_{u} z=-\lap z- \Pd f'(u) \Pp z,
  \]
it is clear that $L_u$ has the required properties: it takes $\Wp$ to $\Wd$ and $\Vp$ to $\Vd$ and the restriction to $v+\Wp$ equals $DF_v$, which is invertible. Moreover, the restriction to  \Vp\ coincides with $-\lap$. This map is no longer a differential operator, due to the integrals needed to compute the projections $P$. But those new terms are innocuous in the finite element formulation --- the sparsity of the underlying matrices is preserved, together with the possibility of standard preconditioning routines.
\begin{figure}[h]
\centering
\mbox{\subfigure{\includegraphics[width=.45\linewidth]{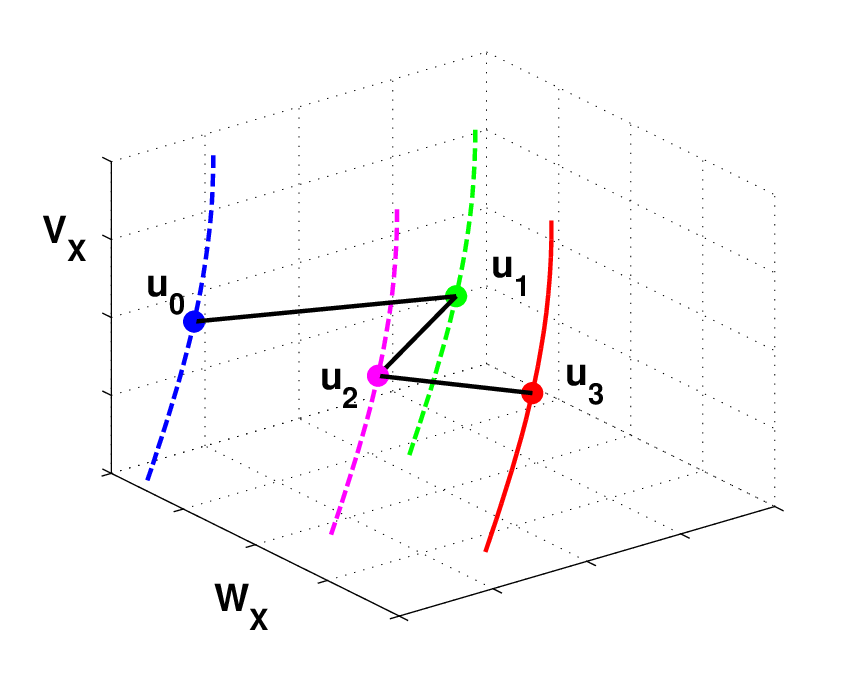}}\quad\quad
\subfigure{\includegraphics[width=.45\linewidth]{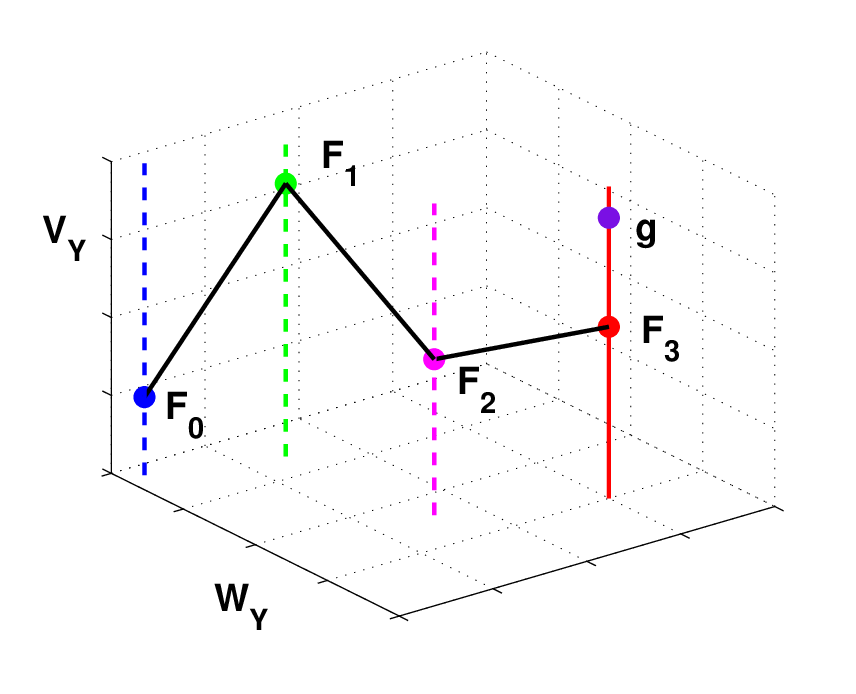}
 }}
\caption{Finding the right fiber}
\label{fig:findfib}
\end{figure}

We search for a point of a horizontal affine subspace $v + \Wp$ which belongs to  $\alpha_g, g\in Y$.The algorithm is straightforward: see Figure~\ref{fig:findfib}. Choose a starting point $u_0$ and
consider its image $F_0$. All would be well if the projections of $F_0$ and $g$ on the horizontal subspace $\Wd$ were equal or at least very near. When this does not happen,
proceed by a continuation method to join both projections. Notice that the algorithm searches for the fiber (i.e., for a point in the fiber) by moving horizontally in the domain.
A direct Newton iteration does not work necessarily: think of finding the (trivial) root of $\arctan(x) = 0$ starting  sufficiently far from the origin.

\section{Moving along the fiber}\label{sec:movfib}

\begin{figure}[h]
\centering
\mbox{\subfigure{\includegraphics[width=.45\linewidth]{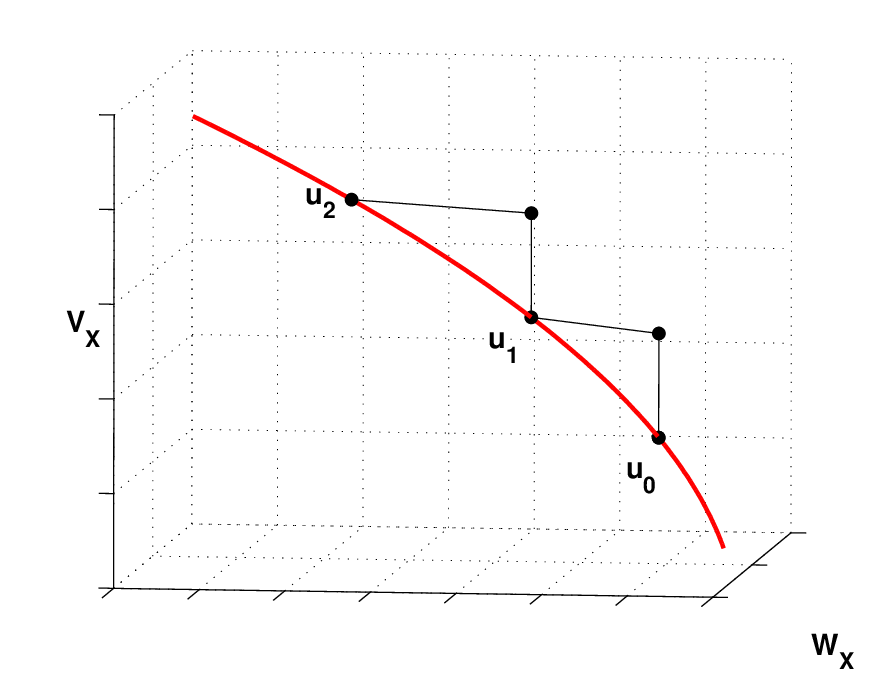}}\quad\quad
\subfigure{\includegraphics[width=.45\linewidth]{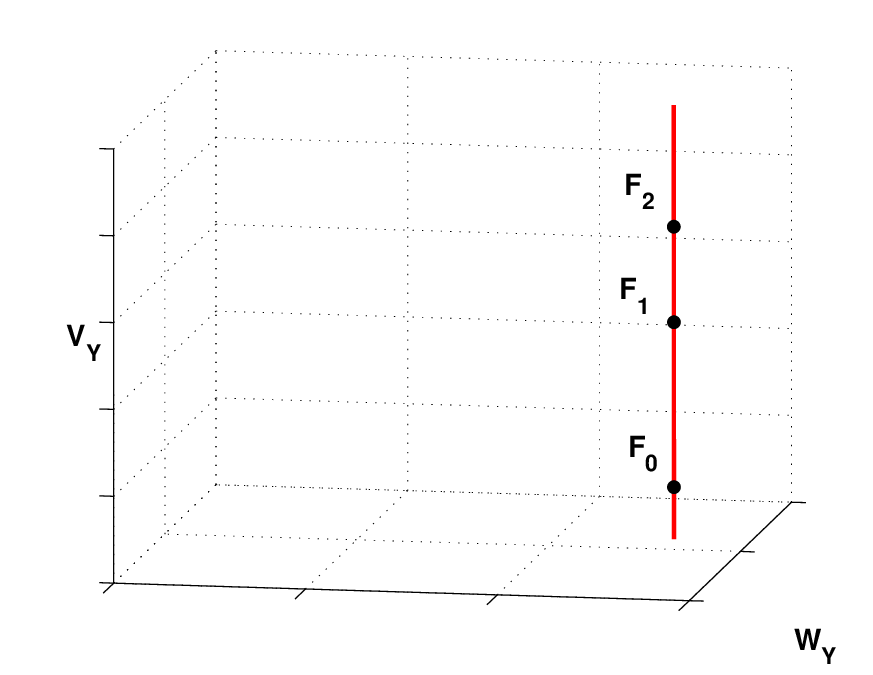} }}
\caption{Mapping a  1-D fiber} \label{fig:mapfib}
\end{figure}

The necessary ingredients for a simple predictor-corrector method to move along a fiber are now available. Say $u \in \alpha_g$ and we want to find another point in $\alpha_g$. Recall that fibers are parameterized by height $v \in \Vp$. Take $u + v$, which is probably not in $\alpha_g$, as a starting point for the algorithm in Section \ref{sec:findfib} to obtain the point of $\alpha_g$ in the same horizontal affine subspace of $u+v$ (see Figure~\ref{fig:mapfib} for two such steps).

We don't know much about the behavior of $F$ restricted to a fiber: the hypotheses on the nonlinearity $f$ are not sufficient to imply properness of $F$, for example. In particular, it is not clear that the restrictions of $F$ to fibers are also proper.

\section{Stability issues}\label{sec:stab}

Proposition \ref{prop:bound} in Section \ref{sec:basics} ensures geometric stability, in the sense that the global Lyapunov-Schmidt decompositions preserve their properties under perturbations. This is convenient when replacing the vertical subspaces spanned by eigenfunctions by their finite element counterparts.

As for the algorithm itself, the identification of the fiber is robust, being a standard continuation method associated to a diffeomorphism between horizontal affine spaces. The numerical analysis along a fiber is a different matter, and the fundamental issue was addressed by Smiley and Chun \cite{smiley:2000}: they showed that the finite element approximations to the restriction of the function $F$ to (compact sets of) the fiber can be made arbitrarily close to the original map in the appropriate Sobolev norm. Here one must proceed with caution: a small metric perturbation may alter the number of solutions, as when switching from $x \mapsto x^2, x \in \RR$ to $x \mapsto x^2 - \epsilon $,  which is a  perturbation of order $\epsilon$ for arbitrary $C^k$ norms. Still, solutions of $F$ which are regular points are stable: they  correspond to nearby solutions of sufficiently good finite element approximations $F^h$.

A different approach might be to interpret the algorithm as a provider of good starting points for Newton's iteration or at least a continuation method. As stated in \cite{plum:2009}, computer assisted arguments require good approximations for the eventual validation of solutions.

\section{Some examples}\label{sec:eg}

For the examples that follow,
$F(u)=-\lap u -f(u) = g,$ with Dirichlet conditions on $\Omega=[0,1]\times[0,2]$. Here, $-\lap_D$ has simple eigenvalues and
$ \lambda_1 = \frac{5}{4} \pi^2 \approx 12.34 ,$ $\lambda_2 = 2 \pi^2 \approx 19.74,$ $\lambda_3 = \frac{17}{4} \pi^2 \approx 41.95.$
\begin{figure}[h]
\centering
\mbox{\subfigure{\includegraphics[width=.45\linewidth]{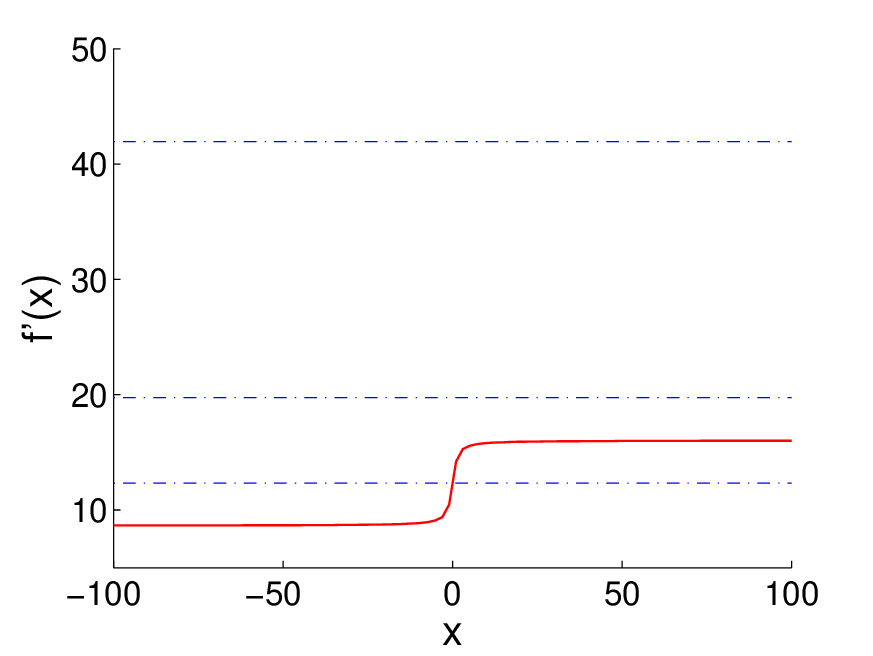}}\quad\quad
\subfigure{\includegraphics[width=.45\linewidth]{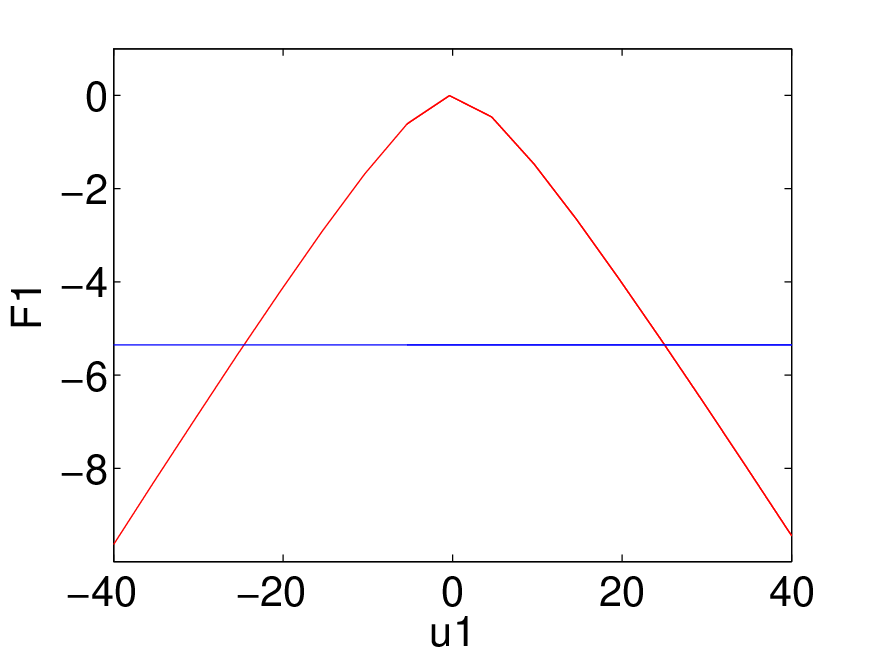} }}
\caption{The derivative of $f$ and the image of $\alpha_g$}  \label{fig:APdffib}
\end{figure}
Denote by $\phi_k^X$ and $\phi_k^Y$ the eigenfunctions of $-\lap_D$ normalized in $X$ and $Y$.

\begin{figure}[h]
\centering
\mbox{\subfigure{\includegraphics[width=.45\linewidth]{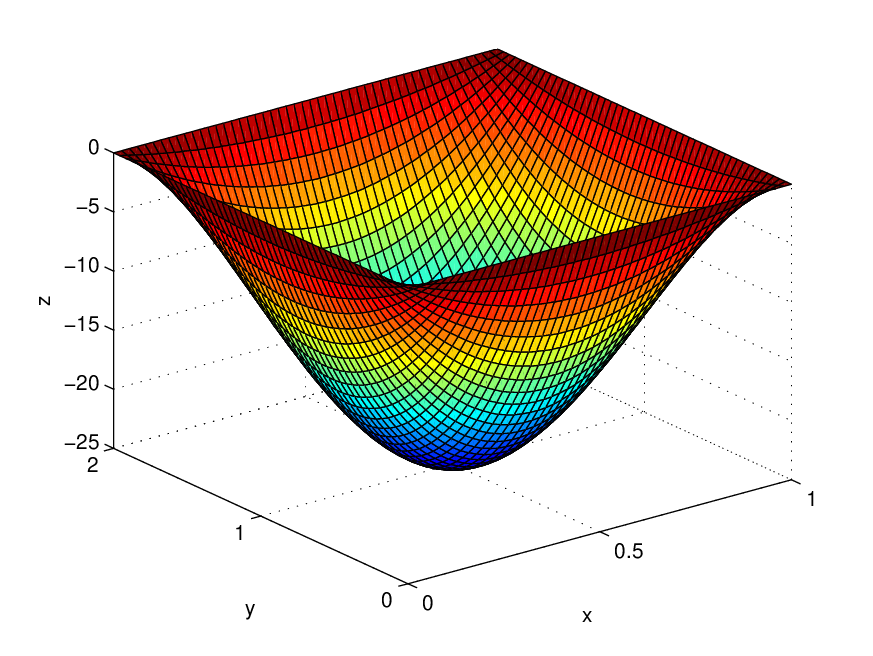}}\quad\quad
\subfigure{\includegraphics[width=.45\linewidth]{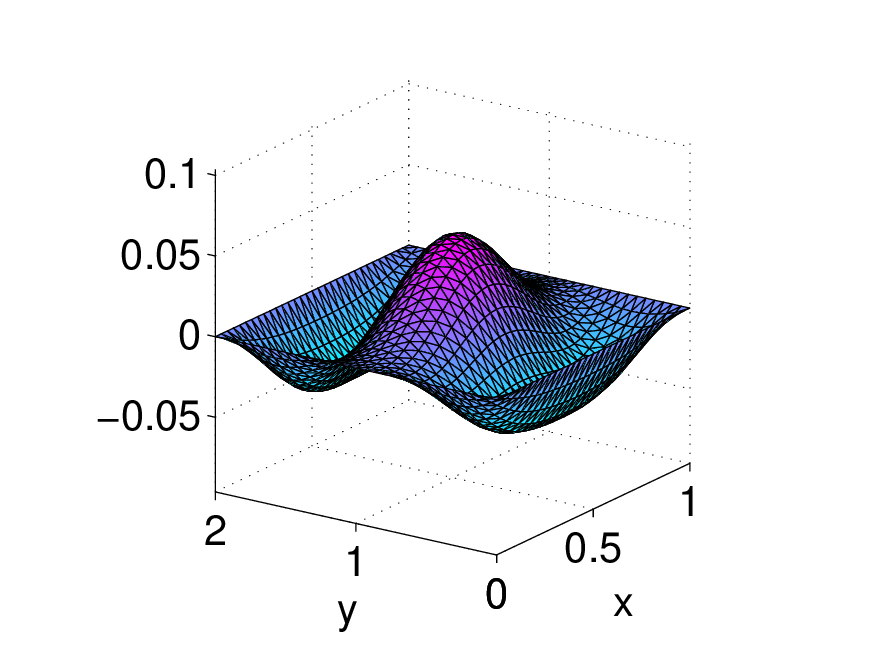} }}
\caption{A right-hand side and a function on its fiber} \label{fig:rhsandfinfib}
\end{figure}

The first example is a nonlinearity $f$ satisfying the hypotheses of the Ambrosetti-Prodi theorem with $f(0)=0$ and derivative $f'(x)=\alpha\arctan(x)+\beta$ with
\[
\Range(f') = \left( \frac{3 \lambda_1 - \lambda_2}{2}\, , \, \frac{\lambda_1 + \lambda_2}{2} \right) >0.
\]
The graph of $f'$ is shown on the left of Figure~\ref{fig:APdffib}.
Here, the index set associated to $f$ is $J=\{1\}$: $\Vp$ and $\Vd$ are spanned by $\phi_{1}^{X}, \phi_{1}^{Y} \ge 0$.
For right-hand side set $g(x) = - 100 x(x-1)y(y-2)$, shown on the left of Figure~\ref{fig:rhsandfinfib}, which has a large negative component along the ground state.

\begin{figure}[h]
\begin{center}
\subfigure{\includegraphics[width=.45\linewidth]{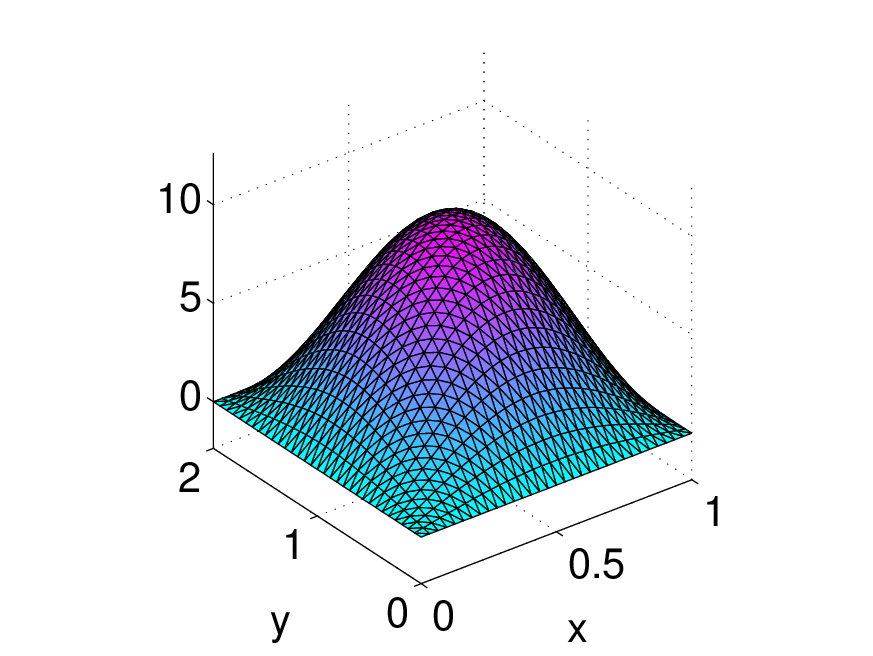}}\quad
\subfigure{\includegraphics[width=.45\linewidth]{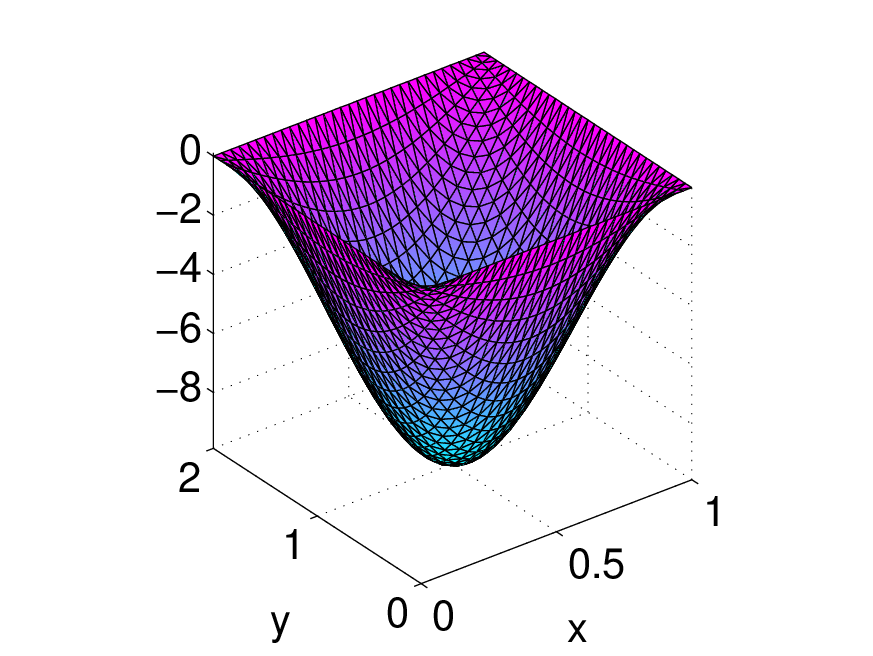} }
\end{center}
\caption{Ambrosetti-Prodi Solutions}
\label{fig:APsols}
\end{figure}

We search for an element of the fiber $\alpha_g$ in the horizontal subspace $\Wp$, starting from $u_0 = 0$, in the notation of Section \ref{sec:findfib}. The result is the function on the right of Figure~\ref{fig:rhsandfinfib}.
Now move along $\alpha_g$, as in Section \ref{sec:movfib}.
The graph on the right of Figure~\ref{fig:APdffib} plots $\iprod[X]{u}{\phi_{1}^{X}}$ (the height of $u\in\alpha_g$) versus $\iprod[Y]{F(u)}{\phi_{1}^{Y}}$ (the height of $F(u)$).
The horizontal line indicates the height of $g$: the solutions of the original PDE correspond to the intersections between the curve and this line.
The two solutions found in this case are presented in Figure~\ref{fig:APsols}.

\begin{figure}[!h]
\centering
\mbox{\subfigure{\includegraphics[width=.4\linewidth]{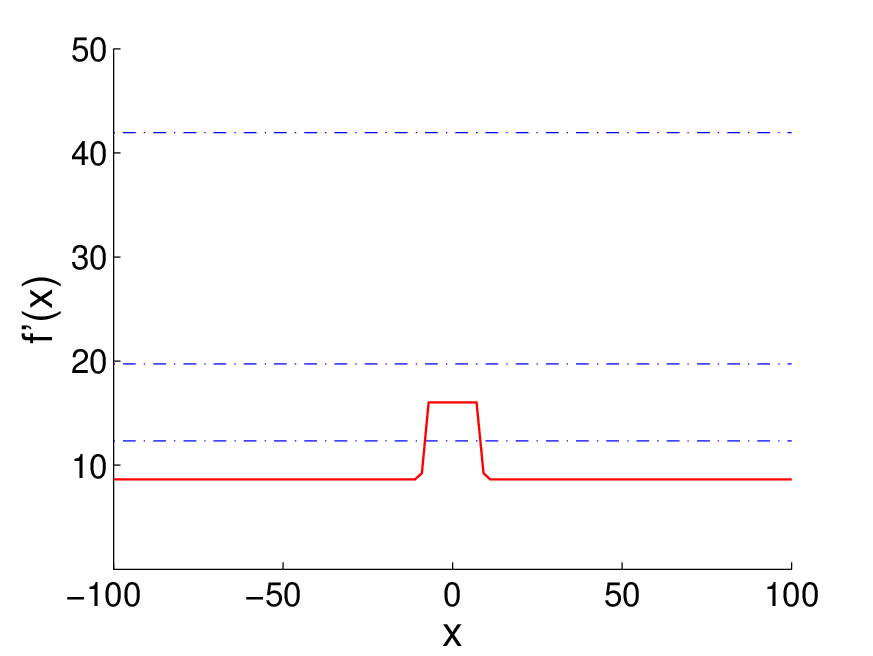}}\quad\quad
\subfigure{\includegraphics[width=.4\linewidth]{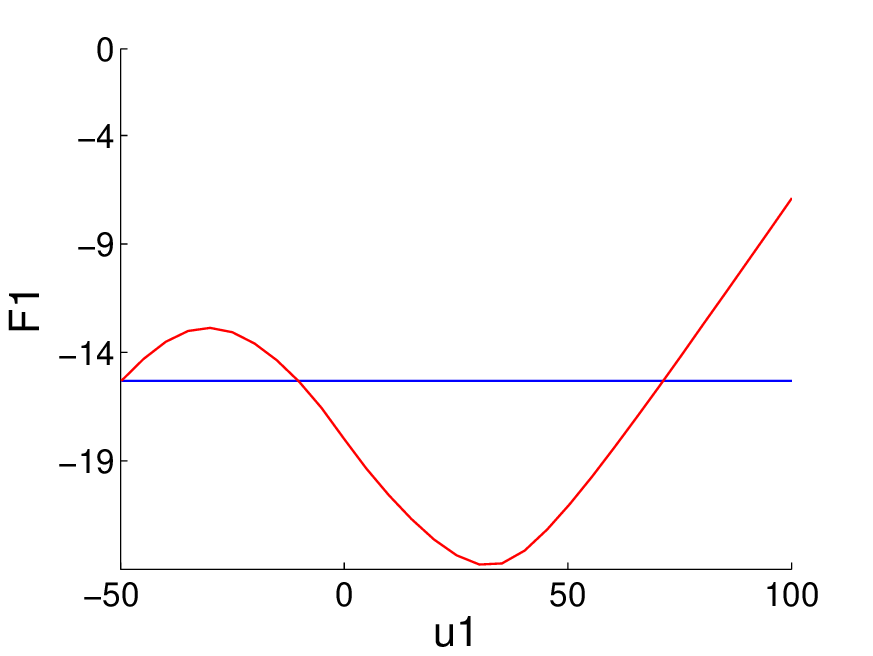} }}
\caption{Non-convex $f$}  \label{fig:NC}
\end{figure}

For the next example, $J=\{1\}$ but $f$ is a nonconvex function whose derivative is depicted on the left of Figure~\ref{fig:NC}.
We consider the fiber through
$u_0(x)=-50\,\phi_1^X(x)+10\,\phi_2^X(x)$, i.e., $\alpha_{F(u_0)}$. According to Figure~\ref{fig:NC} (right),
moving up the fiber yields three distinct solutions.

\begin{figure}[h]
\begin{center}
\includegraphics[width=4cm]{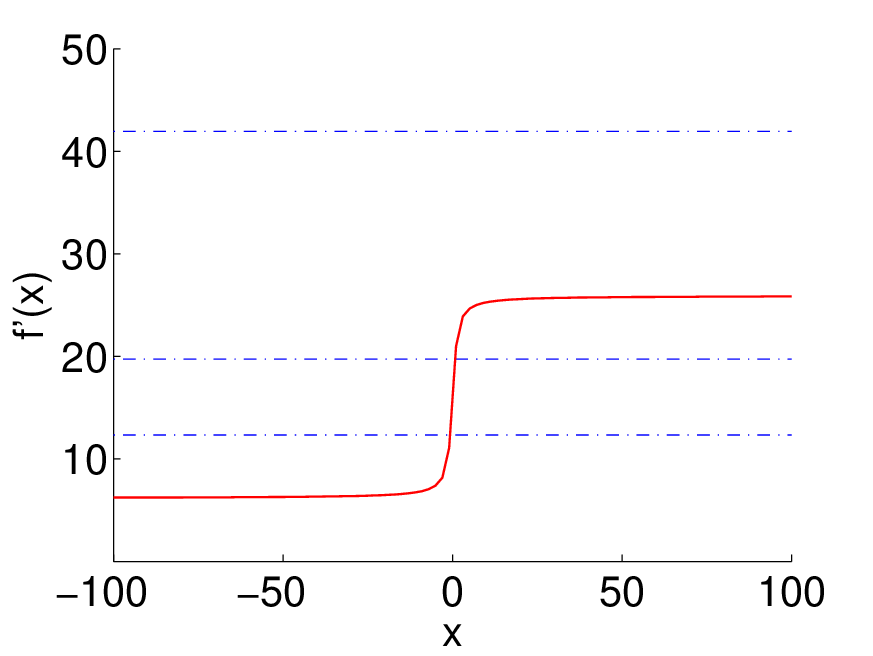}
\end{center}
\caption{The range of $f'$ contains $\lambda_1$ and $\lambda_2$}
\label{fig:APl1l2}
\end{figure}

\begin{figure}[h!]
\centering
\mbox{\subfigure{\includegraphics[width=.4\linewidth]{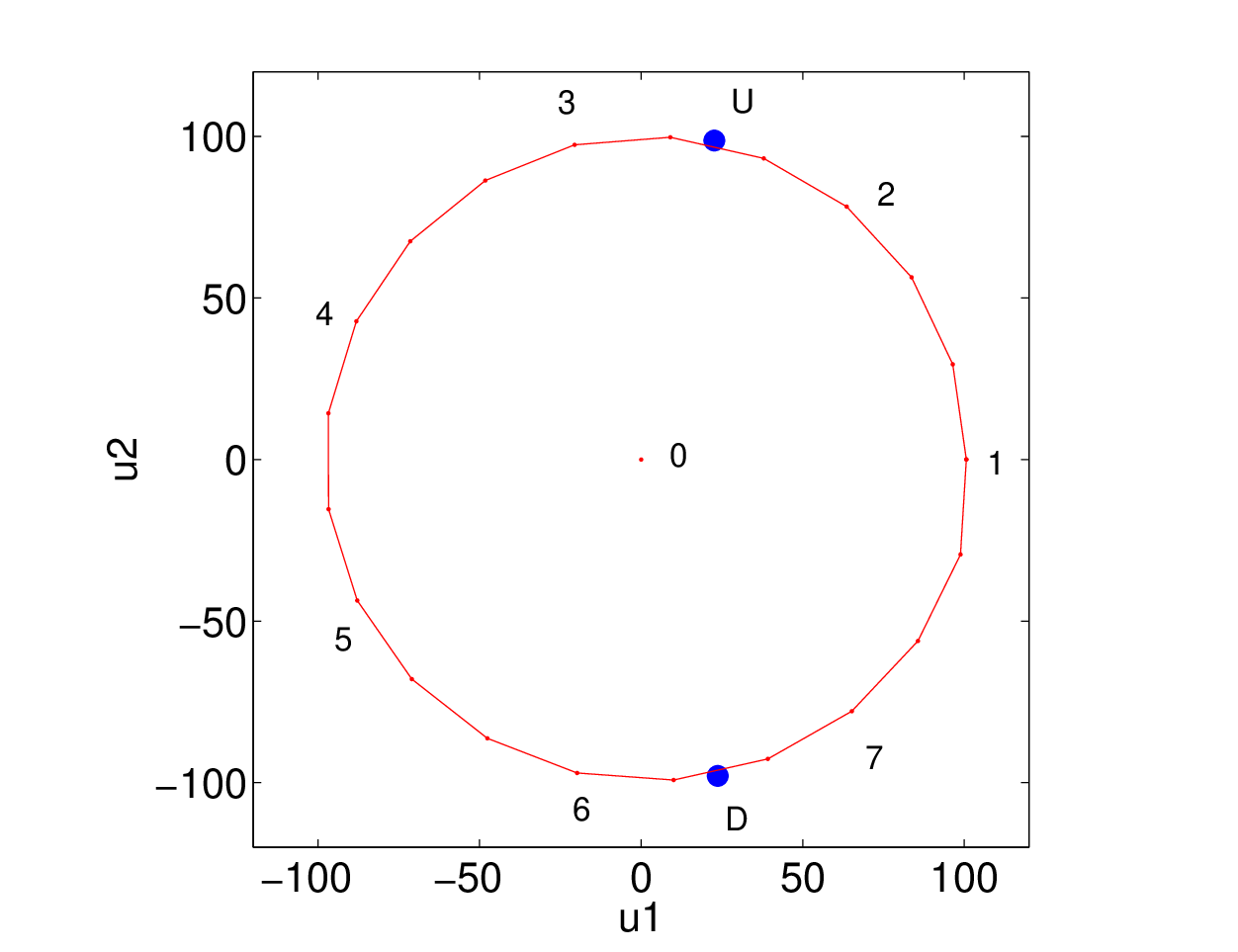}}\quad\quad
\subfigure{\includegraphics[width=.4\linewidth]{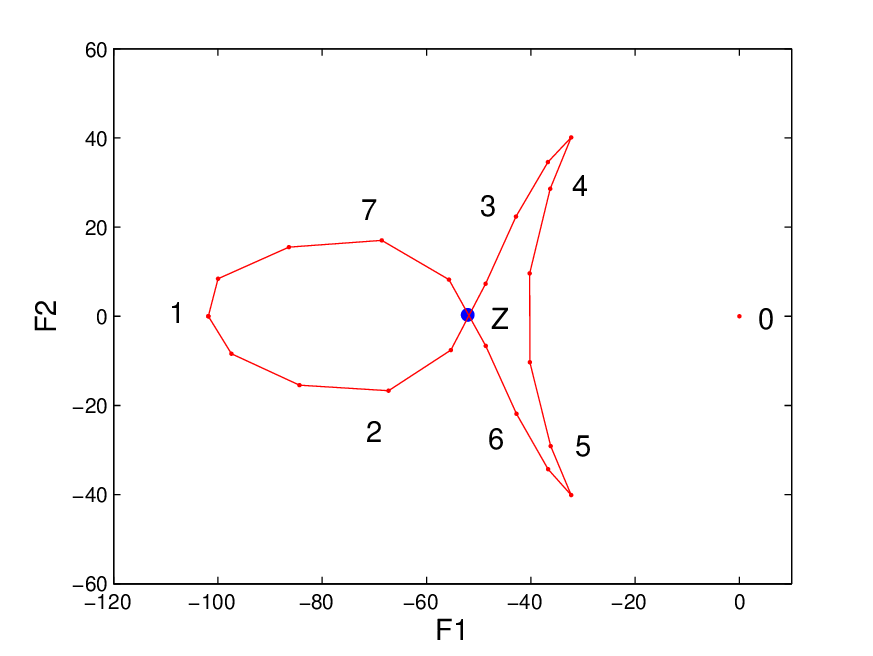}
 }}
\caption{Two preimages on the circle} \label{fig:UD}
\end{figure}

As a concluding example, we take a nonlinearity $f$ for which $J=\{1,2\}$: here the vertical spaces are spanned by the first two eigenfunctions. The function $f$ is of the same form as the first example and its derivative is shown in Figure~\ref{fig:APl1l2}. We study the fiber through the point $u_0 = 0$, which is $\alpha_0$, since $F(0)=0$.
Recall from Section \ref{sec:findfib} that there is exactly one point of $\alpha_0$ for each height, i.e., given a point  $u \in \Vp$, there is a unique point  $\zeta(u) \in \alpha_0$ in the same horizontal affine subspace as $u$.
For a circle $C$ centered at the origin in $\Vp$, $\zeta(C) \in \alpha_0$. The image $F(\zeta(C))$ is shown in the right side of Figure~\ref{fig:UD}: here, we must project $F(\zeta(u))$ along directions $\phi_1^Y$ and $\phi_2^Y$. Clearly, there is a double point $Z$ in $F(\zeta(C))$ and it is not hard to identify in $C$ its two pre-images, $U$ and $D$ marked in the left of Figure~\ref{fig:UD}.

\begin{figure}[h]
\centering
\mbox{\subfigure{\includegraphics[width=.4\linewidth]{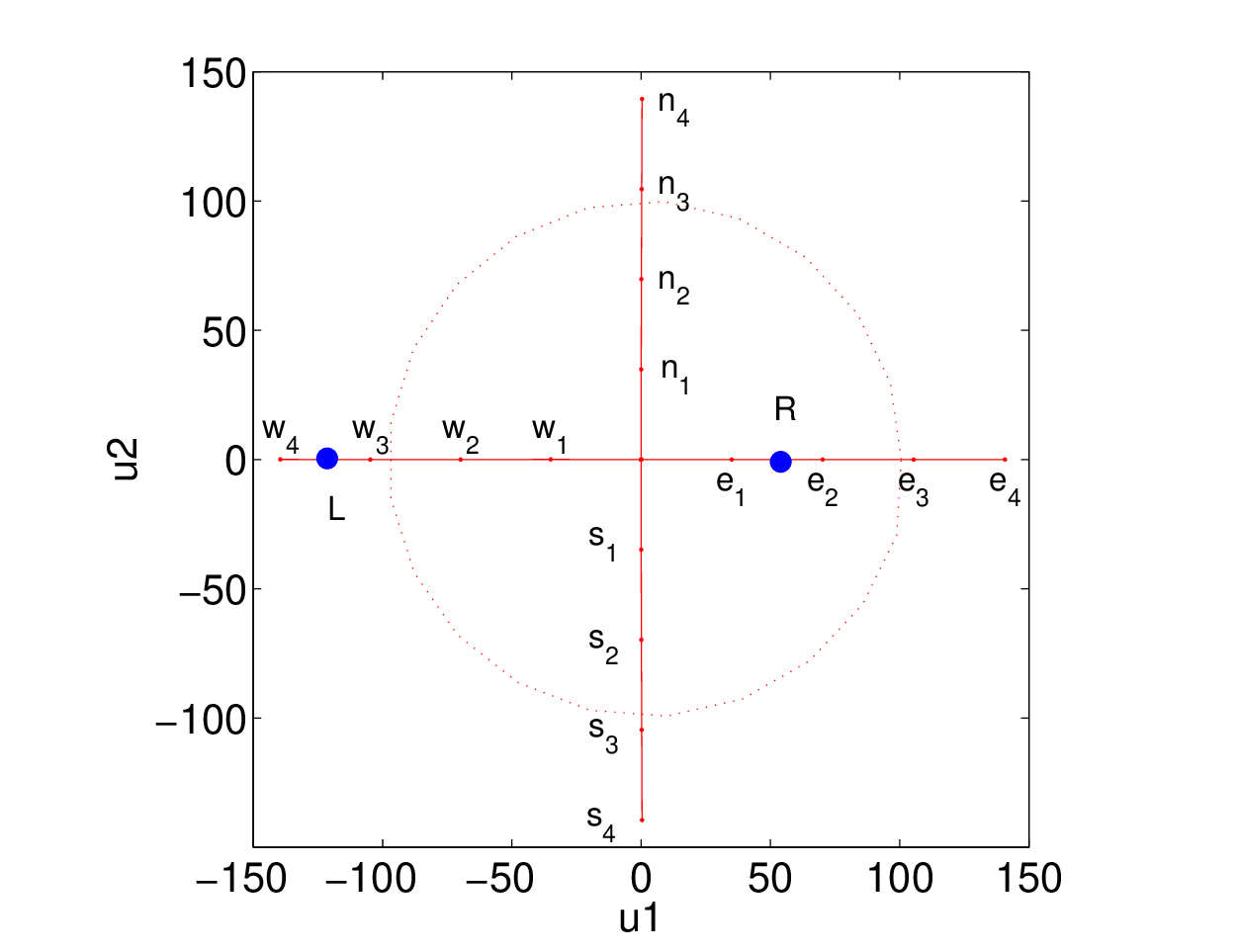}}\quad\quad
\subfigure{\includegraphics[width=.4\linewidth]{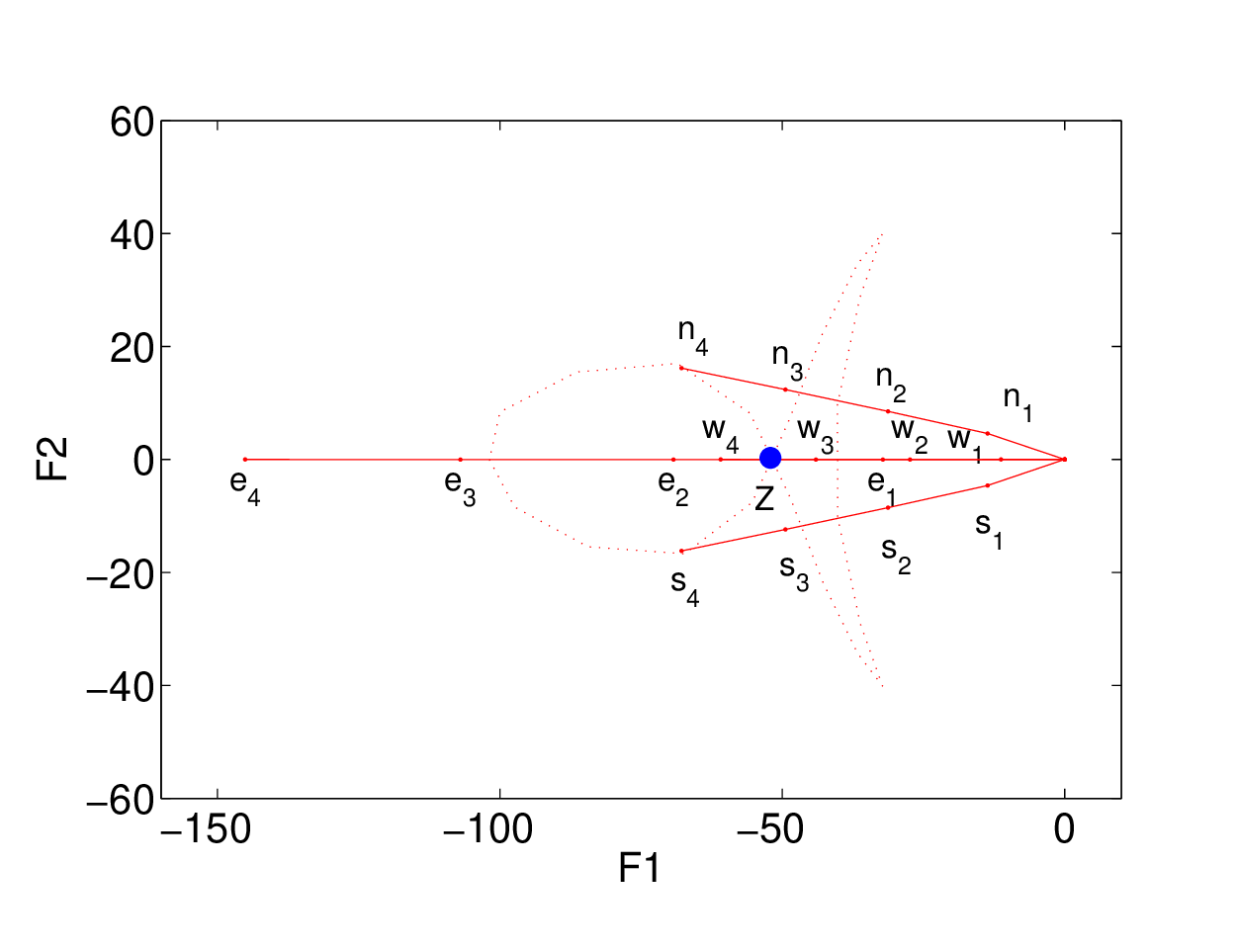}
 }}
\caption{Two preimages along the horizontal axis} \label{fig:LR}
\end{figure}
We now obtain two additional preimages of $Z$ in a rather naive fashion. The images under $F \circ \zeta$ of the four half-axes of $\Vp$ are drawn on the right of Figure~\ref{fig:LR}. It is clear, then, that the horizontal axis contains two preimages $L$ and $R$ of $Z$, which are easily computed.

\begin{figure}[h]
\begin{center}
$%
\begin{array}{cc}
\includegraphics[width=4cm]{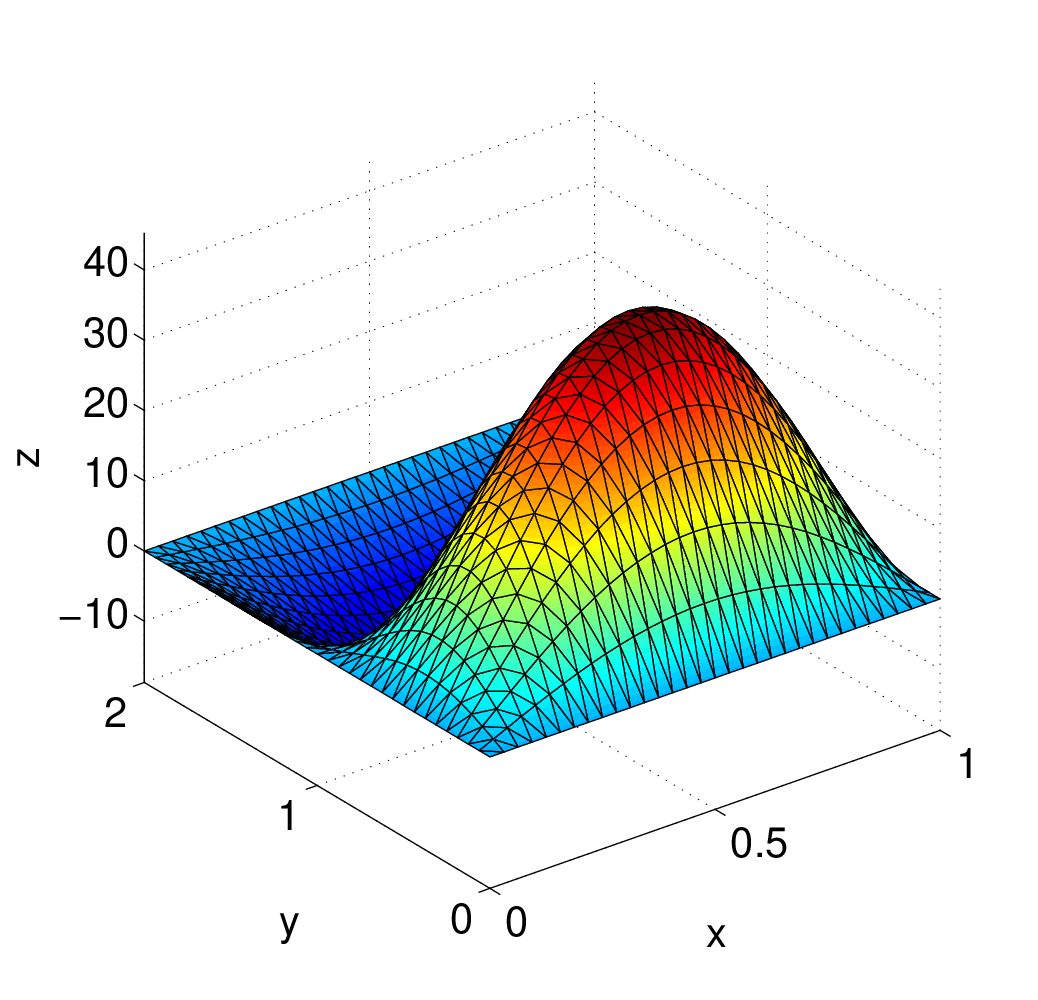} & \includegraphics[width=4cm]{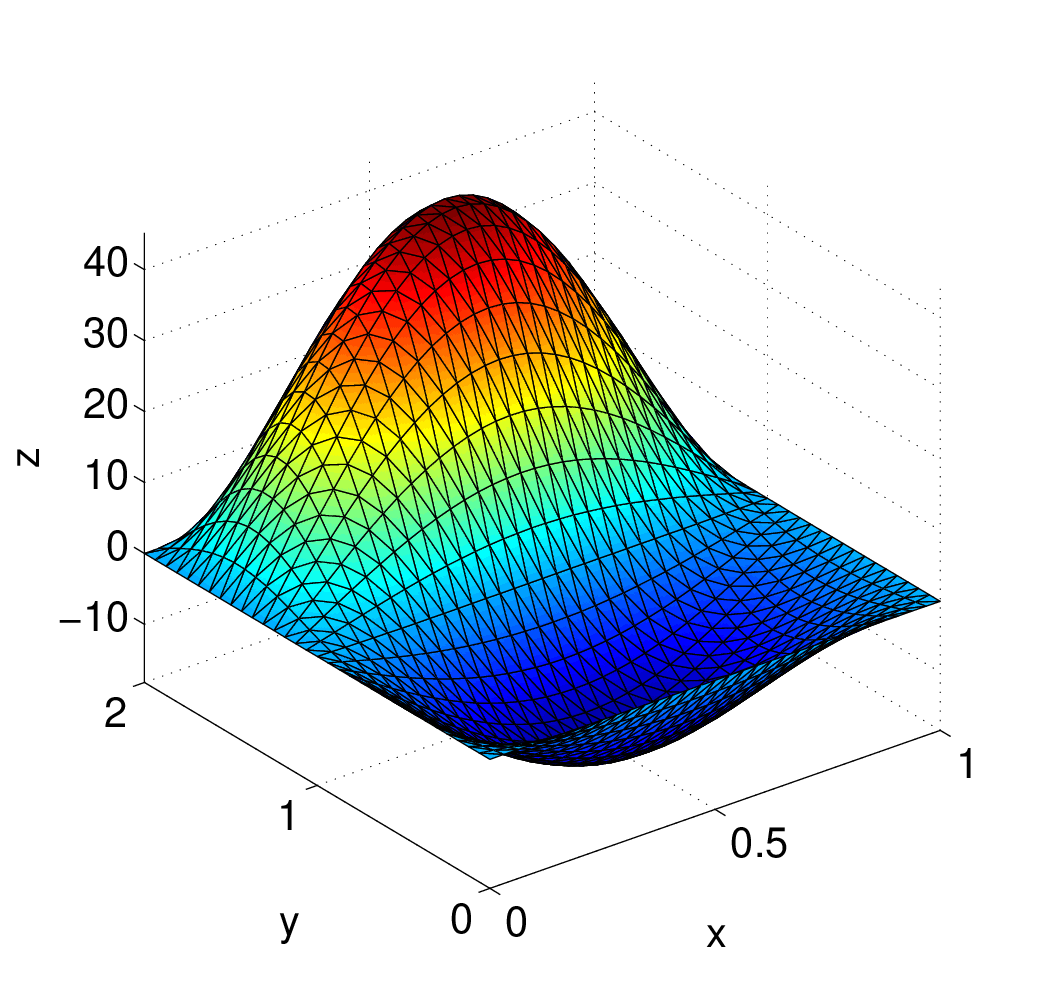} \\
\includegraphics[width=4cm]{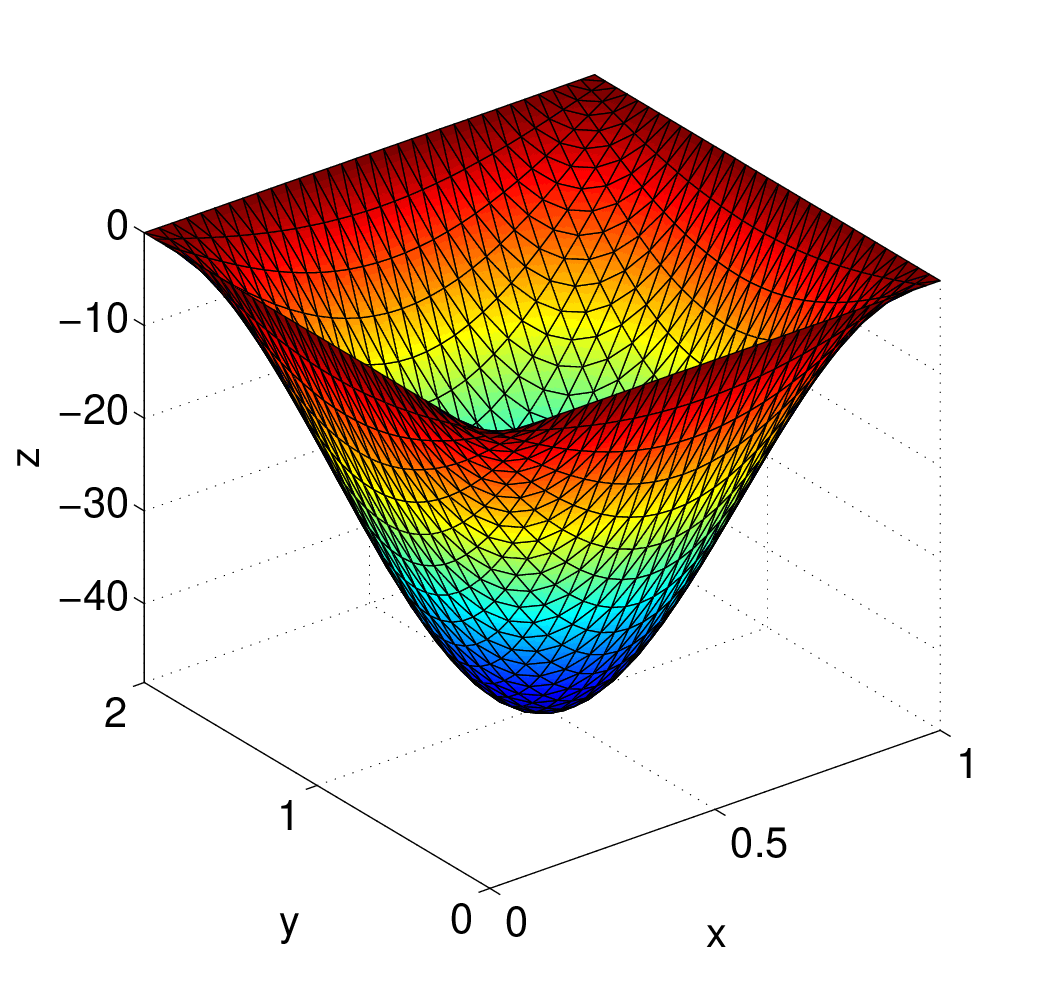} & \includegraphics[width=4cm]{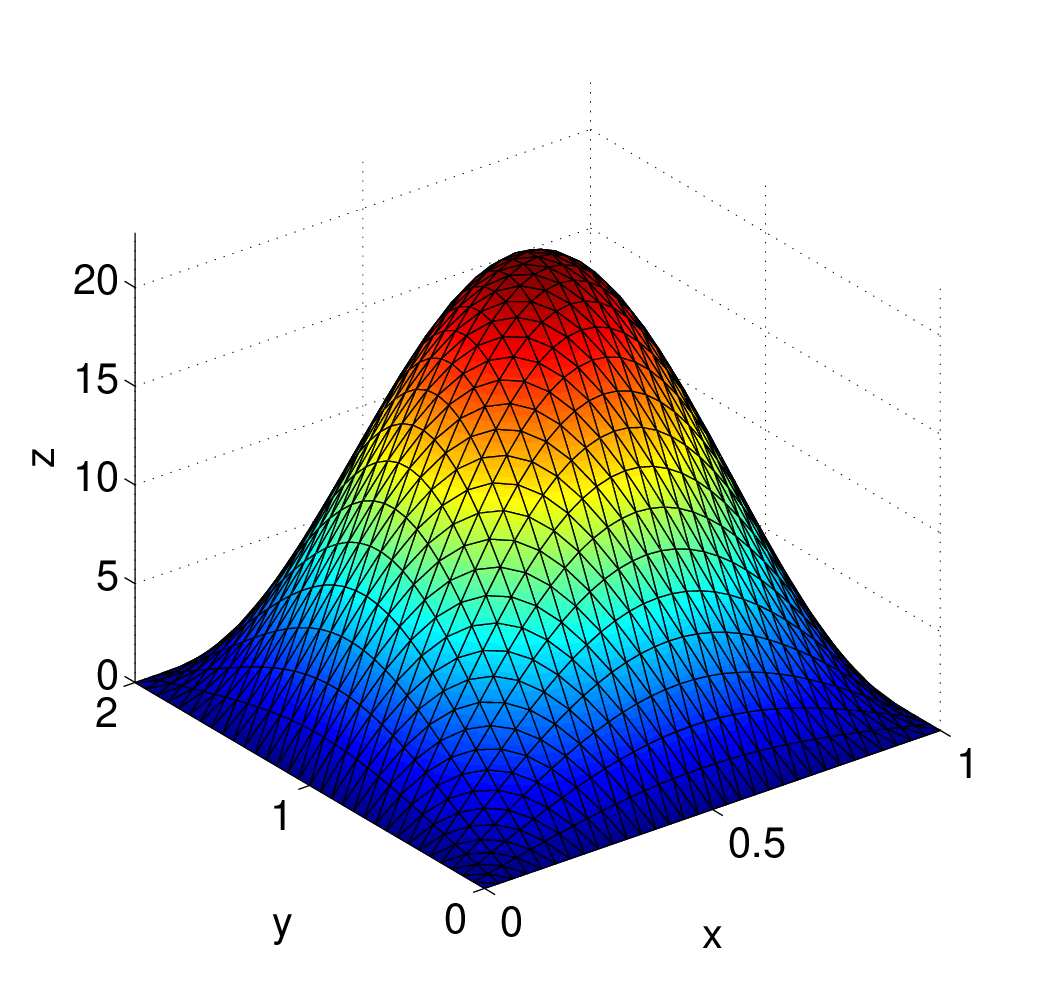}
\end{array}%
$%
\caption{The four solutions}
\label{fig:foursolutions}
\end{center}
\end{figure}

For the sake of completeness, Figure~\ref{fig:foursolutions} displays the four solutions.

\bibliographystyle{ieeetr}

\bibliography{CFL80}		

\end{document}